\documentclass[a4paper,12pt]{amsart}

\usepackage{amssymb,amsbsy,amsmath,amsfonts,amssymb,amscd}
\usepackage{latexsym}
\usepackage{graphics}
\usepackage{color}
\usepackage{fullpage}

\input xy
\xyoption{all}

\newcommand\sC{{\mathcal C}}

\newcommand\sA{{\mathcal A}}
\newcommand\sF{{\mathcal F}}

\newcommand\sH{{\mathcal H}}

\newcommand\la{\lambda}

\newcommand\Ga{\Gamma}
\newcommand\De{\Delta}
\newcommand\ga{\gamma}

\newcommand{\CC}{\ensuremath{\mathbb{C}}}

\newcommand{\ZZ}{\ensuremath{\mathbb{Z}}}

\newcommand{\hol}{\ensuremath{\mathcal{O}}}

\newcommand{\PP}{\ensuremath{\mathbb{P}}}

\newcommand{\ra}{\ensuremath{\rightarrow}}

\def\eea{\end{eqnarray*}}
\def\bea{\begin{eqnarray*}}

\newcommand\dual{\mathrel{\raise3pt\hbox{$\underline{\mathrm{\thinspace d
\thinspace}}$}}}
\newcommand\qe{\ifhmode\unskip\nobreak\fi\quad $\Box$}       

\def\BOX{\hfill\lower.5\baselineskip\hbox{$\Box$}}

\newcommand\R{\Bbb R}

\newtheorem{theo}[equation]{Theorem}
\newtheorem{remarkk}[equation]{Remark}
\newenvironment{rem}{\begin{remarkk}\rm}{\end{remarkk}}

\newtheorem{defin}[equation]{Definition}

\newtheorem{prop}[equation]{Proposition}
\newtheorem{cor}[equation]{Corollary}
\newtheorem{lem}[equation]{Lemma}
\newtheorem{example}[equation]{Example}

\newcommand{\sP}{\ensuremath{\mathcal{P}}}

\def\La{\Lambda}
\def\la{\lambda}
\def\Ga{\Gamma}
\def\ga{\gamma}
\def\R{\rangle}
\def\L{\langle}

\newcommand{\Proof}{{\it Proof. }}
\begin{document}

\title[Cayley forms and self dual varieties]{ Cayley forms and self dual varieties}
\author{ F. Catanese}
\address {Lehrstuhl Mathematik VIII\\
Mathematisches Institut der Universit\"at Bayreuth\\
NW II,  Universit\"atsstr. 30\\
95447 Bayreuth}
\email{fabrizio.catanese@uni-bayreuth.de}

\thanks{The present work was finished in the realm of the DFG
Forschergruppe 790 ``Classification of algebraic
surfaces and compact complex manifolds''.
The first results of this article were announced at  the 1998 Conference in Gargnano, 
and  later at the 2001 Erice
Conference.}

\date{\today}

\maketitle

{\em  This article is  dedicated to  Slava Shokurov
on the occasion of his $60$-th birthday.}

\tableofcontents
\section*{Introduction}

Cayley forms, according to V. Arnold's paradigm by which no mathematical discovery bears the name
of the mathematician who made it first, are nowadays called Chow forms.

A Chow form is a polynomial $F_X$ in the Pl\"ucker coordinates of a Grassmann manifold 
$G (m-n-1, m)$ such that its zero set $$Z =G (m-n-1, m) \cap  \{F=0 \} $$  is the locus of
projective subspaces which intersect a given  projective variety $ X^n_d\subset \PP^m$
(the classical notation $X^n_d$ means that $X$ has  dimension $n$ and degree $d$).

Cayley (\cite{cayley1}, \cite{cayley2}) introduced this concept in the case where $X$ is a curve in $\PP^3$.

His work was later generalized by  Bertini, Chow and van der Waerden (see
\cite{vdw},
\cite{an}, \cite{G-M}, \cite{catchow}, \cite{gkz} for partial accounts), and nowadays, given a
variety  $ X^n_d\subset \PP^m$ as above, one defines
its {\bf Bertini form} $\Phi_X (H_0, \dots, H_n)$ as the minimal polynomial, multihomogeneous of degree
$d$ in each variable $H_i \in  (\PP^m)^{\vee}$ such that 
 $$\Phi_X (H_0, \dots, H_n) = 0 \Leftrightarrow X \cap H_0 \cap \dots H_n \neq \emptyset.$$
This polynomial is very important for applications to vision imaging, since it provides the
`photographic picture' of $X$ for each projection to $\PP^{n+1}$
(if the projection is given by independent  linear forms $(H'_0, \dots, H'_{n+1})  $, the hypersurface image of $X$ is
defined by the polynomial $\Psi$ such that, if we take $H_i = \sum_j a_{ij} H'_j$,
$\Psi (H_0 \wedge \dots \wedge H_n) =  \Phi_X (H_0, \dots, H_n)$).

Moreover, $X$ is completely determined by $\Phi_X $, and there have been several characterizations
of Bertini forms, for instance there is the characterization by Chow and van der Waerden requiring that

1) there exists a polynomial $F$ in the Pl\"ucker coordinates of the Grassmann manifold 
$G (m-n-1, m)$ such that $ \Phi_X (H_0, \dots, H_n) = F (H_0 \wedge \dots \wedge H_n)$:
any such polynomial $F$ is called a Chow form.

2) $ \Phi_X (H_0, \dots, H_n)$ splits as a product of  forms which are linear in $H_n$
in an algebraic extension
of $\CC (H_0, H_1, \dots H_{n-1})$.

Another characterization was given  later in \cite{catchow}, theorem 1.14.

In our opinion the most exciting characterization was given by Green and Morrison (\cite{G-M}), who extended
the result of Cayley, showing that $F$ is a Chow form if and only if certain equations 
of degree 2 or 3 hold identically on the hypersurface $Z =G (m-n-1, m) \cap  \{F=0 \} $.

 The first motivation of this paper was the attempt to see whether the Chow variety
was indeed definable by equations of degree 2 and 3. The impulse for this came from
the beautiful result of Cayley, which we shall now explain  in more detail. 

In this paper a honest Cayley form (respectively: a tangential Cayley form) shall be a polynomial $F$ in the Pl\"ucker coordinates of
$G(1,3)$, whose zero set $Z \subset G(1,3)$ is the set of the lines intersecting a given space curve $C$
(resp.: the lines tangent to a given surface $S$).

$G(1,3)$ is indeed   Klein's quadric in $\PP: = \PP^5$,
defined by 
 $$Q (p) : = p_{01}p_{23} -  p_{02}p_{13}+  p_{03}p_{12}= 0,$$
and this non degenerate quadratic form identifies $\PP$ with its dual space.

 Cayley's equation is 

$$ \frac{1}{2} \{ F, F \} :=   \frac{\partial F}{\partial p_{01}}\frac{\partial F}{\partial p_{23}} - 
\frac{\partial F}{\partial p_{02}} \frac{\partial F}{\partial p_{13}}
+ \frac{\partial F}{\partial p_{03}} \frac{\partial F}{\partial p_{12}} = 0,$$

and Cayley showed that the equation holds on the 3-fold $Z =G (1,3) \cap  \{F=0 \} $
if and only if $F$ is a Cayley form, i.e., either the honest Cayley form of a curve,
or the tangential Cayley form of a surface.

Our  main result  (see theorem \ref{selfdual}) is that this equation 
is equivalent, for a hypersurface $Z \subset
G(1,3)$, to the assertion that $Z$ is selfdual, i.e., $Z$ is equal to its dual variety $Z^{\vee}$.

Examples where a variety and its  dual variety are not hypersurfaces have  for long time been considered,
at least according  to our
knowledge, as sporadic (see \cite{mumford}), and indeed if the variety $X$ is smooth, then Ein (\cite{ein1},
\cite{ein2}) has classified the finite number of cases where $ dim (X) = dim ( X^{\vee})$.

From Ein's classification one can see that there are very few examples where $X$ is smooth
and $X$ and $X^{\vee}$ are projectively equivalent.

Our result says on the other hand that, once we drop the requirement that $X$ be smooth,
 there are countably many families
of self dual varieties, which are not hypersurfaces.

Our second result  expands on a remark made as a footnote to \cite{G-M}, that
a Cayley form (which is not unique) can be changed, by adding a multiple of  Klein's quadric $Q$,
obtaining another Cayley form for which the Cayley equation holds identically on $ Q = G(1,3)$.

We show more precisely (see theorem \ref{Cayley forms}) that there exists a unique representative $F_2$
of the Cayley form such that $ F_2 = F_0 + Q F_1$ with $F_0$ and $F_1$ harmonic,
and such that the Cayley equation for $F_2$ holds identically on the Klein quadric $ Q = G(1,3)$
(i.e., the harmonic projection of the Cayley equation is zero).

This result has as corollary  that the variety of Cayley forms is a projective variety
defined by quadratic equations.

In the same section we also dispose, via elementary examples of curves and surfaces of degree
2 or 3, of too optimistic guesses, that $F_2$ would be just the unique harmonic representative,
or that there exists some representative $F$ such that the  Cayley equation for $F$
is identically zero.

In the final section, we describe (see theorem \ref{honest}) some equations which detect honest Cayley forms
among Cayley forms. These equations appear to be rather simple, however these are again equations which
express that three polynomials vanish identically on the Cayley 3-fold $Z$. The same elementary examples show that  one cannot
alter the Cayley form so that these vanish identically on $Q $, thus showing that the variety of honest Cayley forms is not
a projective variety defined by equations of degree 2 or 3.

 The above results  suggest the question whether the space of generalized Chow forms
(honest and tangential Chow forms) is also defined by quadratic equations. It also suggests the investigation of
the geometric deformations of honest Chow forms to tangential Chow forms. For the time being,
 before finding the solution to this and other questions, we decided to write up this note.

\section{Notation and preliminaries}

Let $V$ be a 4-dimensional vector space over the field $\CC$ (or over an 
algebraically closed field of characteristic $0$), endowed with a volume element ,
i.e., a non zero vector $$ Vol \in \La^4 (V)^{\vee} .$$

The volume element defines a non degenerate symmetric bilinear form 
$$\L, \R \colon   \La^2 (V) \times  \La^2 (V)  \ra \CC: $$ 
$$ \L\omega, \psi \R : = Vol (\omega \wedge \psi). $$

\begin{rem}
The same situation holds for $  \La^m (V) $  when $ dim (V) = 2m$, and $\L, \R $ is symmetric iff $m$ is even,
 skew symmetric iff $m$ is odd.

\end{rem}

In the case where $V = \CC^4$, with canonical basis $e_0, e_1, e_2, e_3$, then we
 have a canonical volume such that 
$Vol (e_0 \wedge e_1\wedge e_2 \wedge e_3) = 1$, and we have, identifying 
$ p \in  \La^2 (V)$ to a skew symmetric $4 \times 4$-matrix
$(p_{ij})$, that one half of the corresponding quadratic form  is just the Pfaffian 
$ Q(p)$ of the  skew symmetric $4 \times 4$-matrix
$$Q (p) : = \frac{1}{2}  \L p,p \R = Pf ((p_{ij})) = p_{01}p_{23} -  p_{02}p_{13}+  p_{03}p_{12}.$$ 

To the symmetric bilinear form $\L, \R$ corresponds the polarity isomorphism
$$
\sP  \colon  \La^2 (V) \ra  \La^2 (V)^{\vee}$$
whose inverse determines a quadratic form on $ \La^2 (V)^{\vee}$, which will be still denoted by $Q$ (this is unambiguous in view of
the polarity isomorphism). 

When $V = \CC^4$, with canonical basis $e_0, e_1, e_2, e_3$, then $ \La^2 (V)$ has canonical basis
$\frac{\partial}{\partial p_{ij}}: = e_i \wedge e_j$, and the quadratic form on $ \La^2 (V)^{\vee}$ yields
the Laplace operator

$$\Delta : =  \frac{\partial}{\partial p_{01}}\frac{\partial}{\partial p_{23}} - \frac{\partial}{\partial p_{02}} \frac{\partial}{\partial p_{13}}
+ \frac{\partial}{\partial p_{03}} \frac{\partial}{\partial p_{12}}.$$

We shall throughout consider polynomial functions $ F (p_{ij})$ on $ \La^2 (V)$, and using the polarity isomorphism
we can define the gradient as the column vector  $ \nabla F$ transpose of the row vector

$$^T \nabla F : = (  \frac{\partial F}{\partial p_{23}},-  \frac{\partial F}{\partial p_{13}}, \frac{\partial F}{\partial p_{12}}, \frac{\partial F}{\partial p_{03}}, - \frac{\partial F}{\partial p_{02}}, \frac{\partial F}{\partial p_{01}}) ,$$
corresponding to the differential $d F$, and  define the Cayley bracket.

\subsection{The Cayley bracket}

\begin{defin}
Let $F (p_{ij}),G (p_{ij})$ be polynomial functions on $ \La^2 (V)$: then their {\bf Cayley bracket} is defined by
the symmetric bilinear form 
$$ \{ F, G \} : = \L \nabla F , \nabla G  \R =  \L d F , d G  \R .$$

The Cayley equation for $F$ is then the differential equation:
$$ \frac{1}{2}  \{ F, F \} = Q (\nabla F) =  \frac{\partial F}{\partial p_{01}}\frac{\partial F}{\partial p_{23}} - 
\frac{\partial F}{\partial p_{02}} \frac{\partial F}{\partial p_{13}}
+ \frac{\partial F}{\partial p_{03}} \frac{\partial F}{\partial p_{12}} = 0.$$ 

\end{defin}

Turning now to geometry, to a homogeneous polynomial $F (p_{ij})$ on $ \La^2 (V)$ corresponds the hypersurface
$$ F : = \{ (p_{ij})  | F (p_{ij}) = 0 \} \subset \PP (\La^2 (V) ) = Proj (\La^2 (V)^{\vee}) \cong \PP^5$$
which we denote by the same symbol $F$.

A particular role plays the hypersurface $Q$, since 
$$ \{ (p_{ij})  | Q (p_{ij}) = 0 \} \subset \PP (\La^2 (V) ) $$
equals the Grassmann manifold
$$ G(1,3) =  \{ p= (p_{ij})  | \exists v, v' \in V, p = v \wedge v' \}$$
parametrizing projective lines $L$  in $\PP(V)  \cong \PP^3$. 

If then $p$ is a point of the hypersurface $F$ (i.e., $ F(p) = 0$), then the tangent hyperplane to $F$ at $p$ is the hyperplane
$$ TF_p : = \{ (\xi_{ij}) | \sum_{ij}  \frac{\partial F}{\partial p_{ij}} \xi_{ij} = 0 \} = \{ \xi | ( dF , \xi ) = 0 \}$$
where $(,)$ denotes the standard duality.

As usual, the non degenerate scalar product $\L , \R$ identifies $TF_p $ to the zero set of the linear form $dF$,
hence to the orthogonal to the gradient  $\nabla F = \sP^{-1} (dF)$.

In particular, if $p \in Q$, then $TQ_p $ is the orthogonal hyperplane $p^{\perp}$ to $p$, since $dQ = \sP (p)$.

In particular , it follows immediately

\begin{lem}
Let $Z$ be the 3-fold in $\PP :=  \PP (\La^2 (V)) $ which is the complete intersection of the Grassmann manifold $Q = G(1,3)$ with
the hypersurface $F$.  
Then the Zariski tangent space to $Z$ at $ p \in Z$ is
$$ TZ_p = p^{\perp} \cap (\nabla F (p))^{\perp}.$$

\end{lem}

We come now to a key formula

\begin{lem}\label{Euler}
Let  $F$ be a homogeneous polynomial of degree $m$ on $\La^2 (V) $. Then  Euler's formula reads out as:

$$   \{ F, Q \} = \L \nabla F,  \nabla Q \R = \sum_{ij}  \frac{\partial F}{\partial p_{ij}}p_{ij} = m F .$$

\end{lem}

\Proof
We have $ \nabla Q (p) = p$, since $ d Q = \sP (p)$,
hence $   \{ F, Q \} = \L \nabla F,  p \R = ( dF , p) = m F$.

\qed

An important consequence is: for $p$ on the hypersurface $F$, one has $ \L \nabla F,  p \R = 0$.

\subsection{Lines in the Grassmannian}

In the sequel we shall denote $\PP(V)$ by $\PP^3$, and by $\PP$ the projective space 
$ \PP (\La^2 (V)) $
 containing the Grassmann manifold $ Q = G(1,3)$ parametrizing lines $ L \in \PP^3$.
We shall use the notation $x,y$ for points in $\PP^3$, and $\pi, \pi '$ for planes in $\PP^3$.

Given $ x \in \PP^3$, $\PP^2_x \subset Q $ is defined as the projective plane in $\PP$, 
$$\PP^2_x : = \{ L  | \ x \in L \} \cong \PP^2,$$ 
and given a plane $\pi \subset 
 \PP^3$, $\PP^2_{\pi}  : = \{ L \  L \subset \pi \} \cong \PP^2$.
 
 Given $x$, $\pi$, one has  $\PP^2_{\pi} \cap \PP^2_x = \emptyset $ unless $ x \in \pi$,
 and in this case one obtains a Schubert line in $ \PP$:
 $$ \Ga (x, \pi) : = \PP^2_{\pi} \cap \PP^2_x  = \{ L  | x \in L \subset \pi \}\ ( \cong \PP^1 \ for \ x \in \pi ).$$
 
 Observe that any line $\Ga \subset Q$ is of this form, and one can find $x, \pi$ as follows.
Let $L, L'$ be two points of $\Ga$, so that the corresponding lines  $L, L' \subset \PP^3$ are not skew
(else $\L L, L' \R \neq  0$): hence $x$ is the intersection point of two corresponding two lines, and
 $\pi \subset \PP^3$ is the plane
spanned by $L, L'$.

 We recover the planes $ \PP^2_x $ and $\PP^2_{\pi}$  starting from $\Ga$ in the following way.
 Intersect $Q$ with
the orthogonal
$\Ga^{\perp}$, and observe that $\Ga
\subset
\Ga^{\perp}$
 is then the vertex of the quadric $Q' : =  Q \cap  \Ga^{\perp} \subset  \Ga^{\perp}  \cong \PP^3$.

 Hence $Q'$ splits as the union of two planes meeting along $\Ga$, which therefore are of
 the form $ \PP^2_x $ for  $ x \in \PP^3$ as above, respectively $\PP^2_{\pi}$ for the above plane $\pi \subset 
\PP^3$.

\subsection{Harmonic polynomials}

Consider the coordinate ring of  $\PP$, namely, the symmetric algebra of 
$ \La^2 (V)^{\vee}$
$$ \sA = \bigoplus_{m \geq 0} \sA_m  : =  \bigoplus_{m \geq 0}  S^m ( \La^2 (V)^{\vee}).$$
Inside $\sA_m $ there is the linear subspace of harmonic polynomials
$$ \sH_m : = \{ F \in \sA_m | \Delta (F) = 0 \} $$ 
where $\Delta$ is, as above,
the Laplace operator

$$\Delta : =  \frac{\partial}{\partial p_{01}}\frac{\partial}{\partial p_{23}} -
 \frac{\partial}{\partial p_{02}} \frac{\partial}{\partial p_{13}}
+ \frac{\partial}{\partial p_{03}} \frac{\partial}{\partial p_{12}}.$$

We recall some basic formulae, which are easy to establish, for homogeneous polynomials $A,B$
(indeed, 3) was proven in lemma \ref{Euler}):

 $$1) \ \ \ \De (AB) = \De (A) B +  A \De (B) +  \L \nabla A, \nabla B\R$$

$$ 2) \ \  \De (Q) = 3 $$

$$3) \ \   \L \nabla A, \nabla Q\R =  deg (A) \cdot A,$$
hence finally
$$1^*) \ \ \ \De (G Q) = (deg (G) + 3)\cdot G + Q \De (G) ,$$
which is the main tool to prove the following

\begin{lem}
There is an isomorphism $\sH_m \cong H^0 (\hol_Q (m)): = W_m$, and moreover 
one has the direct sum decomposition
 $$\sA_m = \bigoplus_{i \geq 0}   Q^i\ \sH_{m-2i} .$$
\end{lem}
\Proof
One  shows the assertion by induction on $m$, using that $\sA_m \cong 
W_m \oplus Q \sA_{m-2}$.

Assume that $G$ is harmonic and let $deg (G) = m - 2i  $; then, by induction on $i$, we easily get:
$$ \ \ \ \De (G Q^i) = i (m + 2 - i)\cdot G \cdot Q^{i-1} .$$
This formula, and the induction assumption shows that the subspaces $Q^i\ \sH_{m-2i}$
build a direct sum inside $\sA_m$, since no harmonic polynomial can belong to the subspace 
$Q \sA_{m-2}$.

Hence there is an injective linear map  $\sH_m \ra  W_m$, and to conclude that
 it is an isomorphism it suffices
(either to show that both spaces have the same dimension, or) to use that both spaces are representations
of $ GL (V)$, and that $W_m$ is irreducible (being the space of sections of a
linearized  line bundle on an homogeneous variety).

\qed

\section{Cayley forms and self dual 3-folds}

\begin{defin}
We shall say that $ F \in H^0 (\hol_{\PP}(m))$ is a Cayley form if the 3-fold $ Z : = Q \cap F = G(1,3) \cap F$
is such that each of its irreducible components $W$ is either 

i) a honest Cayley 3-fold, consisting of the lines $L$ which intersect an irreducible curve $C \subset \PP^3$, 
( $W = \cup_{x \in C} \PP^2_x $) or

ii) a tangential Cayley 3-fold, consisting of the closure of the set of  lines $L$ which are tangent to an irreducible 
non degenerate surface $S \subset \PP^3$ (i.e., $S$ is not a plane) at a smooth point $x \in S$
( $W = \overline{ \cup_{x \in S \setminus Sing(S)} \Ga (x, TS_x) } $).

\end{defin}

\begin{rem}\label{secondderivative}
In the case where $F$ is a honest Cayley form, then $ m = deg(F) = deg(C)$.

If $F$ is a tangential Cayley form associated to a surface $S \subset \PP^3$, then  $ m = deg(F) $
is the intersection number of $ Z : = Q \cap F = G(1,3) \cap F$ with a line $\Ga$ contained in $Q$,
which is then of the form $\Ga(x, \pi)$.

If one denotes by $C'$ the intersection of $S$ with a general plane $\pi$,
one sees therefore that $m$ is the class of the plane curve $C'$.
Thus we have  $$ m =  n (n-1) - \sum_{ y \in Sing(C') } c(y)$$
where $ n = deg (S)$, and $c(y)$ is the Pl\"ucker defect of the singular point $ y \in C'$.

\end{rem}

The following is our first result

\begin{theo}\label{selfdual}

Let  $ F \in H^0 (\hol_{\PP}(m))$, and assume that $ Z : = Q \cap F$ is  reduced. 

Then the following conditions are equivalent:

1) $F$ is a Cayley form,

2) $F$ satisfies the weak Cayley equation $ \{ F, F\} \equiv 0 \ ( mod (Q, F))$,

3)  the 3-fold $ Z : = Q \cap F = G(1,3) \cap F$ is self dual, i.e., $ Z = Z^{\vee}$.

\end{theo}
The structure of the proof runs as follows: first we show  that we can restrict to the case where $Z$ is irreducible,
and  we prove that $ 1) \Rightarrow 2)$;
then we show $2) \Leftrightarrow 3)$, and finally $ 3) \Rightarrow 1)$.

 {\em Proof of theorem \ref{selfdual}, part I.}

Assume that the hypersurface $Z$ is reducible: then we can write $ Z = Z_1 \cup Z_2$
hence, since $Pic (Q) \cong \ZZ$, changing $F$ modulo $Q$, we may assume $ F = F_1 F_2$,
with $ F_1,  F_2$ relatively prime.

Then $$ \{ F, F\} = \L  d F , dF \R =  \L  F_1 d F_2 +  F_2 d F_1,  F_1 d F_2 +  F_2 d F_1 \R =$$
$$ F_1 ^2  \{ F_2, F_2\} + 2 F_1 F_2 \{ F_1, F_2\} + F_2 ^2  \{ F_1, F_1\}.$$

Hence  $F_1$ and $F_2$  satisfy 2) if and only if $F$ does.
Therefore we may restrict ourselves to show the theorem in the case where $Z$ is irreducible.

$ 1) \Rightarrow 2)$:

{\bf Case i) where $F$ is a honest Cayley form of an irreducible curve $C$.}

Let $ L \in Z$: then there is $ x \in C$ such that $ L \in \PP^2_x \subset Z$, hence $F$ vanishes
on $\PP^2_x$. Take now coordinates on $\PP^3$ such that $ x = e_0$, hence
$\PP^2_x = \{p |  p_{12} = p_{13} =p_{23} =0 \}$, whence $ \nabla F (L)  $
has components which satisfy 
$$ \frac{\partial F}{\partial p_{0i}} (L) = 0 , \ i=1,2,3 \Rightarrow  \{ F, F\} (L) = 0.$$
Thus $ \{ F, F\} $ vanishes on $Z$, equivalently the weak Cayley equation 2) holds.

{\bf Case ii) where $F$ is a tangential  Cayley form.}

Let $ L \in Z$ be general: then there is $ x \in S$ which  is a smooth point and is such that $L$ is tangent to $S$ at $x$.
 Take now coordinates on $\PP^3$ such that $ x = e_0$, $ L = e_0 \wedge e_1$, and the tangent
 space $TS_x$ is the plane $\{ x | x_3= 0 \}$.
 
 There exists a local parametrization of $S$ with
 $$ x = ( 1, u,v, \phi (u,v))$$
 where $\phi$ has order at least two at the origin $ u=v=0 $.
 
 Then a local parametrization for the variety of tangent lines is given by the wedge product of
 the two (row) vectors:
 $$ ( 1, u,v, \phi (u,v) )$$
 $$  ( 0, 1,\la , \phi_u (u,v) + \la \phi_v (u,v) )$$
 
  hence the lines are parametrized by $(u,v, \la)$, $L$ corresponds to the origin in this system of coordinates,
  and we have
$$  p_{01} = 1, p_{02} = \la, p_{03} = \phi_u (u,v) + \la \phi_v (u,v) ,  p_{12} = u  \la - v,$$
$$  p_{13} = u (  \phi_u (u,v) + \la \phi_v (u,v)) - \phi (u,v).$$

Notice that, since $  p_{01} = 1$, $  p_{23} =   p_{02} p_{13} - p_{03} p_{12}$ on $Q$
and looking at the Taylor development of the function
$$ F ( p (u,v, \la)) =  \frac{\partial F}{\partial p_{02}}(L) \la + \frac{\partial F}{\partial p_{03}}(L)  \phi_u (u,v) -  \frac{\partial F}{\partial p_{12}}(L) v +  {\rm terms \ of \ order \geq 2},$$
which is identically zero, we obtain that, at the point $L$, $\frac{\partial F}{\partial p_{02}}$ vanishes, 
and $\frac{\partial F}{\partial p_{03}}$ vanishes too unless  $\phi_{uu} (0,0): =  \frac{\partial ^2 \phi}{\partial u ^2} (0,0) = 0$.

Moreover $\frac{\partial F}{\partial p_{01}}(L)$ vanishes by  Euler's formula.

The conclusion is that $ \{ F, F\} (L) = 0$ unless the tangent line $L$ is a zero of the II fundamental form of $S$
(a so called asymptotic direction).
But since the surface is non degenerate, for general $L$ we have that $L$  is not a zero of the II fundamental form of $S$.

Hence  $ \{ F, F\} $ vanishes on $Z$, equivalently the weak Cayley equation 2) holds.

\qed

The above calculation in local coordinates shows that, if $L$ is a smooth point of $Z$, then the tangent space $TZ_L$
is the subspace $ \{ p | p_{13}= p_{23} = 0 \}$, which contains the $\PP^2_x$ of lines passing through $x$.

It also shows the following

\begin{prop}\label{secondderivative}

 If the line $L$ is not an asymptotic direction at $x \in S$, then the second derivative of $F$ does not identically vanish
on $\PP^2_x$.

\end{prop}

\Proof
$\PP^2_x$ is the subspace  $ \{ p | p_{12}= p_{13} = p_{23} = 0 \}$,
and we are claiming that the second fundamental form of $Z$ does not vanish on it.

Intersecting $Z$ with this subspace we obtain the subvariety defined by
$$  v = \la u, \ \ \ \ u (  \phi_u (u,v) + \la \phi_v (u,v)) = \phi (u,v) \Leftrightarrow $$  
$$ \Leftrightarrow  v = \la u,  \ \ \ \ u   \phi_u (u,\la u) +   \la u  \phi_v (u,\la u)) - \phi (u,\la u)=0. $$

All we have to show is that at the origin the function  
$$  u   \phi_u (u,\la u) +   \la u  \phi_v (u,\la u)) - \phi (u,\la u) $$ has a quadratic term which is not  identically zero.

But this quadratic term equals the one of 
$$ u   \phi_u (u,\la u)  - \phi (u,\la u) .$$

Letting $ \phi (u,v) = a u^2 + b uv + c v^2 \ ( \ mod \ (u,v)^3 )$, we obtain $$ u ( 2 a u ) - a u^2 = a u^2 \equiv 0,$$
hence  $0 =  2 a = \phi_{uu} (0,0) $, contradicting our assumption.

\qed

 {\em Proof of theorem \ref{selfdual}, part II.}

$ 2) \Leftrightarrow 3)$:

2) just says that, for $ L \in Z$, $ Q (\nabla F (L) ) = 0$: this means that $\nabla F (L) $ is a point in $Q$.

However, since $$  \L \nabla F (L) , \nabla F (L)  \R = 0 ,  \L L, L   \R  = 0,  \L \nabla F (L), L  \R =0, $$
where the last equality is nothing else than the Euler formula (Lemma \ref{Euler}), we see that
2) is equivalent to saying that the line $\Ga_L :  L * \nabla F (L)$  joining $L$ and $\nabla F (L)$ is fully contained in the Grassmannian $Q$.

Observe now that, identifying $\PP$ with its dual space via the polarity $\sP$, the line
$\Ga_L := L * \nabla F (L)$ is dual to the pencil of tangent hyperplanes to $Z$ at $L$: since
 $TZ_L  = L^{\perp} \cap \nabla F (L)^{\perp}  $.

We have therefore shown the following

{\bf Claim: 2) holds $\Leftrightarrow$ we have the inclusion of the dual variety of $Z$ in $Q$: $$ Z^{\vee} \subset Q.$$ } 

We conclude the proof of this step via part 2) of the following lemma.

\begin{lem}
Assume that $ Z \subset Q$. Then

 1)  $ Z \subset Z^{\vee}$.
 
 2) $ Z^{\vee} \subset Q \Leftrightarrow  Z =  Z^{\vee}.$

\end{lem}

{\em Proof of the Lemma.}

1):  assume that $ L \in Z$ is a smooth point: then $TZ_L \subset TQ_L =  L^{\perp}$. Hence $ L \in Z^{\vee}$.

2): $ Z^{\vee} \subset Q   $ implies, by  1), that $ Z^{\vee} \subset (Z^{\vee})^{\vee} = Z$, where the last equality is the
biduality theorem. Again by  1)  $ Z \subset Z^{\vee}$, hence $ Z^{\vee} \subset Q   $ implies $Z =  Z^{\vee}$,
while the converse is obvious.

\qed

\qed

The following proposition explains the geometrical background for the last step of proof of Theorem \ref{selfdual}.
It involves the concept of Segre dual curve,  that we need to recall (see \cite{piene}: however, for the reader 's benefit,
we give an elementary  proof).

\begin{defin}
Let $C$ be a non degenerate curve in $\PP^n$, which means that, if $\ga(t)$ is a parametrization
of $C$, then for general $t$ the $n$ vectors $ \ga(t), \ga'(t)  , \dots \ga^{(n-1)}(t) $
are linearly independent.

Then the Segre dual curve $C^* \subset (\PP^n)^{\vee}$ is the curve of osculating $(n-1)$-dimensional spaces,
so that  $C^* $ is parametrized by  $$\ga^* (t): = \ga(t)\wedge  \ga '(t)  \wedge \dots  \wedge \ga^{(n-1)}(t) .$$

More generally, the k-th associated curve $C[k]$  is the curve of osculating $(k)$-dimensional spaces,
a curve in the Grassmann manifold $ G(k,n)$, parametrized by
$$\ga [ k ] (t): = \ga(t)\wedge  \ga '(t)  \wedge \dots  \wedge \ga^{(k)}(t) .$$

\end{defin}

\begin{lem}\label{Segre}
If $C$ is a non degenerate curve in $\PP^n$, then

a) $(C^*)^* = C$

b)   for each value of the parameter $t$, $\ga ^*[n-1-k] (t) $ is the  annullator subspace of $ \ga [k] (t) $ 

c) $C^{\vee} $ is the tangential developable hypersurface of $C^*$.
\end{lem}

\Proof
Observe that a) is the special case of the more general statement b), obtained taking $k=0$.

In order to prove b), we use the method of moving frames. Namely, we let $A(t)$ be the matrix
with columns the $n+1$ vectors  $$ \ga(t), \ga'(t)  , \dots \ga^{(n-1)}(t), \ga^{(n)}(t) .$$
$A(t)$ determines a flag in $\CC^{n+1}$, and we may also take a unitary matrix  $U(t)$ determining the same flag.

Then the `dual flag',  given by the annullators of these subspaces in the dual space  $\CC^{n+1}$, 
corresponds to the matrices  $B(t), V(t)$ where one takes the respective  dual bases in the opposite order.

One considers as usual the Cartan matrix $C(t)$, the skew symmetric matrix defined by 
$$ U^{\cdot} (t): = \frac{d U(t)}{dt} = C(t) U(t).$$

We have that $ ^T V(t) U (t) \equiv J$, where $J$ is the antiidentity matrix ; whence, taking the derivative of both sides,
$$  ^T V(t) U ^{\cdot} (t)  +^T V ^{\cdot}(t) U (t) = 0 \Rightarrow  \ ^T V(t) C(t)  +^T V ^{\cdot}(t) = 0 \Rightarrow $$
$$ \Rightarrow V ^{\cdot}(t) = C(t) V(t) .$$
This formula shows that the dual flag is the osculating flag of the curve $\ga^*(t)$.

One can also avoid the use of the complex numbers, and work with the moving frame $A(t)$, defining the companion matrix
$M(t)$ such that $ A^{\cdot} (t) = M(t) A(t),$ and the proof follows similarly.

To prove the last statement, observe that 
$$C^{\vee} = \{ H | \exists x \in C, TC_x \subset H \} = \{ H | H  \in Ann  \ga [1] (x) \} = $$
$$  \{ H | H  \in   \ga ^*[n-2]  (x) \}  =  \{ H | \exists x \in C, H  \in  {\rm Linear \ span} ( \ga ^*(x), \dots ,  \ga ^{*(n-2)}  (x)) \} .$$

\qed

\begin{prop}

Consider the (involutory) polarity isomorphism identifying $\PP$ with its dual space, which geometrically
corresponds to the mapping associating to a line $L \subset \PP^3$ the pencil of planes containing
it (a line in $(\PP^3)^{\vee}$).

It sends the tangential Cayley 3-fold of a surface $S$ to the tangential Cayley 3-fold of the dual variety  
$S ^{\vee}$ when the latter is a surface $S$, else to the honest Cayley 3-fold of the dual variety  
$S ^{\vee}$ when the latter is a curve.

It sends the honest Cayley 3-fold of a curve $C$ to the  tangential Cayley 3-fold of the dual variety  
$C ^{\vee}$, which is the tangential developable surface of the Segre dual curve $C^*$. 

\end{prop}
\Proof
We use the standard notation by which the projectively dual subspace of a projective subspace $L \subset \PP^n$,
i.e., the projective subspace corresponding to the annullator, is denoted by $L^*$.

Now, if $L$ is a tangent line to the surface $S$ at a point $x$, then $ x \in L \subset TS_x$, hence, defining $H : = TS_x$, 
we have  $ H^*  \in L^* \subset x^*$, thus $L^*$ is tangent to $S^{\vee}$, which settles the proof in the case
where $S^{\vee}$ is a surface (in view of biduality).

Again by biduality, it suffices to consider the honest Cayley 3-fold of a curve $C \subset \PP^3$.
It consists of the lines $L$  intersecting the curve $C$ in a point $x$; then the dual subspace $L^*$
satisfies  $ H^*  \in L^* \subset x^*$, whenever the plane $H$ contains $L$. We choose $H$
to also contain $TC_x$, so that $H^* \in C ^{\vee}$, and $L^*$ is tangent to $C ^{\vee}$ at $H^*$.

Conversely, if $L^*$ is tangent to $C ^{\vee}$ at $H^*$, then there is $x$ such that $ H^*  \in L^* \subset x^*$,
and $x \in L$.

\qed

 {\em Proof of theorem \ref{selfdual}, part III.}

$ 3) \Rightarrow 1)$:

For each smooth point $L \in Z$, the line $\Ga_L : =  (L * \nabla F(L))$  corresponds to the pencil of tangent hyperplanes to $Z$ in $L$,
hence it is contained in  $Z^{\vee} = Z$.

Being a line in the Grassmannian, it determines a point $x \in \PP^3$ and a plane $\pi \subset \PP^3$ such that
$\Ga_L = ( L * \nabla F(L)) = \Ga (x, \pi)$.

Hence, we get a rational map of $Z$ onto a correspondence 
$$\Sigma \subset \PP^3 \times (\PP^3)^{\vee}: = \overline{ \{ (x, \pi) | \exists L  \in  Z \setminus Sing(Z),  \ s. \  t.   \  \Ga_L = \Ga (x, \pi) \}}.$$

\begin{lem}
$\Sigma$ has dimension 2 and is a duality correspondence with respect to the two projections.
\end{lem}

\Proof
For each point $ L \in Z$ , we have the line  $\Ga_L : =  (L * \nabla F(L)) = \Ga (x, \pi)$ which is contained in $Z$. Assume that there is
another line $\Ga '$ passing through $L$, contained in $Z$ and different from $\Ga_L$. Then $\Ga '$ is contained in $TZ_L = \Ga_L^{\perp}$.
Hence the plane $\Pi$ spanned by  $\Ga_L$ and by $\Ga '$ is contained in $TZ_L$ and we have then $\Pi \subset Q$, since 
 $\Ga ' \subset TZ_L = \Ga_L^{\perp}$. 
 
 Since $\Ga (x, \pi) =  \Ga_L \subset \Pi \subset Q$, it follows that either
 
 1] $\Pi = \PP^2_x$, or
 
 2] $\Pi = \PP^2_{\pi}.$
 
 We separate our analysis according to different cases:
 
 i) for general $L \in Z$, there is only  a finite number of  lines passing through $L$ and  contained in $Z$.
 
 ii) for general $L \in Z$ there is an infinite  number of  lines   contained in $Z$ and passing through $L$.
 
Condition ii) implies,  by the above consideration, that one of the following  holds:
 
 [1]: for general $L \in Z$, $L \in  \PP^2_x \subset Z $
 
  [2]: for general $L \in Z$, $L \in  \PP^2_{\pi} \subset Z $.
  
Therefore, if ii) holds true, then necessarily $Z$ is a honest Cayley 3-fold, or a dual honest Cayley 3-fold.

 Consider now the tangential correspondence $W$ for $Z' := Z \setminus Sing(Z)$:
 $$ W : = \{ (L_1, L_2) \in Z' \times Z' | \ TZ_{L_1} \subset L_2^{\perp}  \} =  \{ (L_1, L_2) \in Z \times Z |\  L_2 \in \Ga_{L_1}   \}.$$
 
Since $ dim (W) = 4$, and $Z$ has dimension 3, the general fibre $Y : = W_{L_2}$ of the second projection is irreducible
of dimension 1. And, for each $L_1 \in Y$, $ L_2 \in \Ga_{L_1} $. Since i) holds and $Y$ is irreducible, it follows
that all the lines  $\Ga_{L_1} $ are equal, and the fibre $Y$ equals $\Ga_{L_1} $. In particular, the tangent space to $Z$
is constant along $\Ga_{L_1} $. We also obtain  that  the map onto $\Sigma$
is constant over $\Ga_{L_1} $, hence $\Sigma$ is a surface. 

Moreover since ii) does not hold, the two projections of $\Sigma$ yield two surfaces, $S \subset \PP^3$, $S' \subset ( \PP^3)^{\vee}$.

There remains to show that $S$ and $S'$ are dual to each other.
Now, for each general point $ x \in S$, $x$ is the image of a line $\Ga (x, \pi) \subset Z$. If we show that the lines
$ L \in \Ga$ are tangent to $S$ then this proves that $\pi = \cup_{L \in \Ga} L$ is tangent to $S$ in $x$, hence $S'$ is 
dual to $S$. 

This assertion is proven in the forthcoming Lemma.

\qed

\begin{lem}
Let $ f \colon Z \setminus Sing (Z)  \ra S$ be  the above morphism, such that $ f (L) = x$, where $x$ is the intersection point
of the lines $ L, \nabla \subset \PP^3$, $ \nabla : = \nabla F (L)$. 

Then $ \PP^2_x \subset TZ_L$, and, if $Df$ is of maximal rank at $L$, then $ Df (\PP^2_x) = L$.

\end{lem}
\Proof
Letting as usual $\Ga$ be the line joining $L$ with $\nabla$,
we know that $ TZ_L = \Ga^{\perp} $, that $\Ga \subset  \Ga^{\perp} $, $\Ga \subset Z \subset Q$.

Then  $ TZ_L \cap Q = \PP^2_x \cup \PP^2_{\pi}$, where $\pi$ is the plane spanned by the lines  $ L, \nabla \subset \PP^3$.

View now $  L$ and $\nabla $ as 4x4 skew symmetric matrices, so that $x$ is the solution of the system
$$ L x = 0 ,  \nabla x = 0. $$
Consider a tangent vector to $L$ with direction $ L' \subset \PP^2_x $: then, if we work as usual with the ring
$\CC[\epsilon ]/ ( \epsilon^2)$, we obtain the equation 
$$ ( L + \epsilon L')  ( x + \epsilon x')  = 0 ,  (\nabla + \epsilon \nabla ') ( x + \epsilon x') = 0 $$
for the first order variation of $f$ along the tangent direction $L'$.

Hence we obtain
$$   L'  x + L x'    = 0 ,   \nabla '  x + \nabla x' = 0 \Rightarrow  L x'   = 0,$$
since $ L'  x = 0 $.

The conclusion is that $ x' = Df (L') $ lies in the line $L$. On the other hand $Df$ has maximal rank (=2), and $\Ga$ lies in the
kernel,  hence $ Df$ satisfies $ Df (\PP^2_x) = L$.

\qed

\qed

\begin{rem}
The Cayley 3-folds $Z$ considered above are all singular. In fact Ein (\cite{ein1}) classified the smooth projective
varieties $X$ such that $ dim (X) = dim (X^{\vee})$ (he actually forgot to explicitly mention the assumption of smoothness,
but this is clearly used, see cor. 1.4 of \cite{ein1}).

\end{rem}

\begin{rem}
Igor Dolgachev pointed out another characterization of Cayley forms in terms of singular loci of line complexes
(see \cite{Jessop}, page 308, \cite{Dolgachev}, page 534).

Since it is related to the previous discussion, we give a brief account in our terminology. A line complex is a subvariety
$Z \subset Q = G(1,3)$. 

We denote by $\La = \PP (U)$ the projectivization of the tautological subbundle on the Grassmannian $G(1,3)$.
Hence  $$\La = \{ (x,L) | x \in L \} \subset \PP^3 \times G(1,3).$$
Denote by  $$\La_Z : = \{ (x,L) | x \in L\in Z \} \subset \PP^3 \times Z,$$
the restriction of the bundle to $Z$, and denote by $f$ the projection on $\PP^3$.

While $\La$ is the fibre bundle $ \PP (T_{\PP^3})$, with fibre over $x \in \PP^3$ equal to $\PP^2_x$,
the same does not occur for  $\La_Z$.

The {\bf singular locus} of the line complex is defined to be the critical set $\sC$ of $ f \colon \La_Z  \ra \PP^3$,
while the {\bf focal locus} is by definition $\sF : = f (\sC)$, the set of critical values of $f$.

Therefore the singular locus equals the closure of  the set of pairs $(x,L)$, $L$ being a smooth point of $Z$,
 where the fibre of $f$ is not smooth of the right codimension;
i.e., such that $  \PP^2_x \cap Z $ is not a transversal intersection at $L$.

In the case where $ dim (Z) = 3$, this means that 
$$ \PP^2_x \subset TZ_L = L^{\perp} \cap \nabla^{\perp}  \Leftrightarrow L, \nabla \in \PP^2_x \Rightarrow \nabla \in Q.$$

In particular, $\sC \subset \La_{Z \cap \{ \{F, F\} = 0\}}$. Conversely, proceeding as in the first two lines of the proof of Lemma 15,
one sees that, if $\nabla \in Q$,  then $ TZ_L \cap Q = \PP^2_x \cup \PP^2_{\pi}$,
thus  $\sC$ projects birationally onto $  Z  \cap \{ \{F, F\} = 0\}.$ 

The intepretation pointed out by Dolgachev is therefore that $Z$ is a Cayley 3-fold if and only if
it equals the projection of its  singular locus. 
\end{rem}

\section{Quadratic equations for the variety of Cayley forms}

A Cayley 3-fold is the divisor $Z$ on the Grassmann manifold $ Q = G(1,3)$ of a section $\zeta \in H^0 (Q, \hol_Q (m))$.

A Cayley form $F$ is a homogeneous polynomial of degree $m$, $F  \in H^0 (\PP, \hol_{\PP} (m))$ such that the restriction
of $F$ to the quadric $Q$ is precisely $\zeta$. Hence we may change a given Cayley form $F$ by adding a multiple 
of $Q$ to it, trying to see whether one could obtain a Cayley form satisfying the strong Cayley equation 
$\{ F,F \} \equiv 0$. We shall show that this cannot be achieved, but at least  (as stated in \cite{G-M})
one can obtain $\{ F,F \} \equiv 0 ( mod\  Q).$

We show indeed a more precise result, which has as consequence that  Cayley 3-folds are parametrized by
a projective variety which is the intersection of quadrics.

\begin{prop}\label{modQ}
Assume that $F$ is homogeneous of degree $m$ and satisfies the weak Cayley equation 
$$  \{ F,F \} \equiv 0 \ ( mod \  (F, Q)). $$

Then there exists another Cayley form $F_2$, defining the same Cayley 3-fold $Z$ as $F$, such that 
$$  \{ F_2,F_2 \} \equiv 0 \ ( mod \   Q). $$

Moreover, $F_2$ is unique $ mod \ (Q^2)$.

\end{prop}

\Proof
We seek for $F_2 = F + QG$ and
 calculate (using the formula $ \{ F,Q \} = m F$) 
$$  \{ F_2,F_2 \} =  \{ F + QG ,F + QG \} =$$
$$  \{ F ,F  \} + 2 Q   \{ F, G \}  + 2 m GF + 2 (m-2) Q G^2 + Q^2 \{ G, G \} + 2 G^2 Q.$$

Hence, if 
$$  \{ F ,F  \} = A Q + B F ,$$
it suffices to take $ G = \frac{-1}{2m} B$, and the solution $G$ is unique modulo $Q$, hence 
$F_2$ is unique modulo $Q^2$.

\qed

We reach then as an important consequence the following

\begin{theo}\label{Cayley forms}
The variety $\sC_m$ of Cayley 3-folds $Z \in \PP (H^0 (\hol_Q (m))$ is isomorphic  
to the  subvariety  $\ \sC'_m \subset \PP (\sH_m \oplus Q \sH_{m-2})$ defined by quadratic equations:
$$\sC'_m   : = \{ Z = Q \cap F  |  F \in (\sH_m \oplus Q \sH_{m-2}) , h_{2m-2}( \{ F, F\}) = 0 \} ,$$

 (here $ h_m : \sA_m \ra \sH_m $ is  the harmonic projector).

\end{theo}

\begin{rem}
Let $F = F_0 + Q F_1\in \sH_m \oplus Q \sH_{m-2}$: then the equation 
$h_{2m-2}( \{ F, F\}) = 0$ can be rewritten as 
$$ h_{2m-2} (\{ F_0, F_0\}  + 2 m F_0 F_1 ) = 0. $$ 

\end{rem}

\subsection{ The easiest examples}

Let $F$ be a Cayley form, so that there are polynomials $A, B$ such that
$  \{ F ,F  \} = A Q + B F .$

Take then $F_2 = F + QG$ as above, where   $ G = \frac{-1}{2m} B + C Q$ is as above.

In the special case where $ deg (F) \leq 3$, then we have the unicity of $F_2$, since $ G = \frac{-1}{2m} B $
by degree considerations ( $deg (C) < 0 $). 

Moreover,  $A,B$ are both unique.

We have then 
$$  \{ F_2,F_2 \} = 
A Q  -  \frac{1}{m} Q   \{ F, B \}  +  \frac{m-1}{2m^2} Q B^2 + \frac{1}{2^2m^2} Q^2 \{ B, B \} .$$

Let us start by considering the case $ deg (F) = 2$.

\begin{cor}\label{examples}
In the case of a smooth quadric surface $S_2 \subset \PP^3$ there is no tangential Cayley form
$F$ satisfying the strong Cayley equation $$  \{ F ,F  \} \equiv 0 .$$
The unique Cayley form $F_2$ such that $  \{ F_2,F_2 \} \equiv 0 \ ( mod \    Q) $ 
is  harmonic.

Even worse  occurs for the honest Cayley forms of two skew lines, or of the twisted cubic curve:
there is no tangential Cayley form
$F$ satisfying the strong Cayley equation $$  \{ F ,F  \} \equiv 0 ,$$
moreover the  unique Cayley form $F_2$ such that $  \{ F_2,F_2 \} \equiv 0 \ ( mod \    Q) $ 
is  not harmonic.

In the case of a smooth plane conic curve, instead,  the harmonic representative satisfies the strong Cayley equation.

\end{cor}

\Proof
  Take the tangential Cayley form of
the quadric surface $$ S = \{ x | x_0 x_1 - x_2 x_3 = 0 \}.$$

A direct calculation shows that a Cayley form is given by
$$ F := (p_{01} + p_{23})^2 + 4 p_{03} p_{12},$$
and that 
$$  \{ F ,F  \} =  8 F .$$

We obtain ( since then $A=0$, $ G = - 2$)
$$  \{ F_2,F_2 \} =    8  Q.$$
Hence 
$F_2 = F - 2 Q$, and  $ \De (F_2)  = \De (F - 2 Q) = 6 - 6 = 0.$

Actually, as pointed out by Dolgachev, if  we are starting from a  quadric surface which  is diagonal with equation
$\sum_i a_ix_i^2= 0$, the corresponding form  is $F = \sum_{ij} a_ia_jp_{ij}^2 $, which is directly seen to be harmonic,
moreover one has $$  \{ F ,F  \} =  4 a_0 a_1 a_2 a_3  Q .$$

In the case of the honest Cayley form of a conic, a Cayley form is easily calculated as
$$ F := p_{02} ^2 + 4 p_{01} p_{12},$$
which is easily seen to be harmonic and to satisfy the  strong Cayley equation.

If we instead take two skew lines, then a Cayley form is
$$ F :=  p_{01} p_{23},$$
satisfying $ \De F = 1$, $  \{ F ,F  \} =   2 F .$ Hence its harmonic representative is $ F - \frac{1}{3} Q$,
while $F_2 = F - \frac{1}{2} Q$, which satisfies 
$  \{ F_2,F_2 \} =    \frac{1}{2}   Q.$

In the case of the twisted
cubic curve, a Cayley form $F$ is obtained as the  determinant of the following symmetric matrix:

$$
\left( \begin{matrix}
p_{01}& p_{02} & p_{03}\\
p_{02} & p_{12}+ p_{03} & p_{13}\\
p_{03}  & p_{13} & p_{23}\\
\end{matrix}\right)
$$

An easy calculation shows that
$$ \De (F) =   p_{12}+ p_{03}.$$

Hence, if $ F = F_0 + Q F_1$ is the harmonic decomposition of $F$, then $4  F_1 = \De (F) =   p_{12}+ p_{03}.$

We skip the rest of the explicit calculations, using a limiting argument: the twisted cubic admits as a limit a chain of 
3 lines, with Cayley form 
$$ F :=  p_{01} p_{02} p_{23},$$
we get:
$$  \{ F ,F  \} =   2 F p_{02} \Rightarrow F_2 = F - \frac{1}{3} p_{02}  Q$$
hence
$$  \{ F_2,F_2 \} =   -  \frac{2}{3}   Q \{ F, p_{02}  \} + \frac{4}{9}  Q p_{02} ^2 = 0 +  \frac{4}{9}  Q p_{02} ^2 .$$
Finally, observing that (since $F$ has degree $3$) $F_2$ is here unique:
$$ \De (F_2 ) =  \De (  F - \frac{1}{3} p_{02}  Q)  =   p_{12}+ p_{03} - \frac{4}{3}  p_{02} . $$
\qed 
\noindent

\section{Equations for honest Cayley forms}

In the previous sections we have shown that the space of Cayley forms is a projective variety
defined by quadratic equations.

Our geometrical explanation shows also that in this variety there are three sets:
1) the closed set of honest Cayley forms (the Cayley forms of some curve $C$  in $\PP^3$)

2) the closed set of dual honest Cayley forms  (the Cayley forms of the developable surface 
$S$ dual
to some curve $C'$ in $(\PP^3)^{\vee}$)

3) the open set of tangential and dual tangential Cayley forms (here $S$, $S^{\vee}$
are  both surfaces).

We are therefore looking for equations which define the smaller closed sets, in particular the first one.

A simple way to obtain such equations is to observe that, while for honest Cayley forms the Cayley
3-fold 
$Z$ contains the $\PP^2_x$ determined by $L$, for a tangential Cayley 3-fold this space is
contained in 
$ TZ_L$ (indeed $ TZ_L \cap Q = \PP^2_x \cup \PP^2_{\pi}$) but, according to proposition  
\ref{secondderivative}, the second derivative  of $F$ does not vanish on $\PP^2_x$ for general $L \in Z$.

Therefore we want that for $ L \in Z = \{ L | Q(L) = F(L) = 0 \}$, the quadratic form
$D^2 F(L) (p, p) $ associated to the Hessian matrix of $F$ vanishes identically
on $$\PP^2_x = \{ p |   p \wedge x (L) = 0\}.$$

To have explicit equations use the following elementary 

\begin{lem}
Let $L, L' \in Q$ be two coplanar lines in $\PP^3$ such that the plane $\pi$ spanned
by them does not contain the point $e_0$.
Then, letting $x$ be the intersection point of the two lines, the plane
 $\PP^2_x$ has as basis $L, L'$ and $ L'' = e_o \wedge x$.

Writing $ L'' = \sum_{i=1}^3 y_ i e_0 \wedge e_i = e_0 \wedge y$, we obtain that
the Pl\"ucker coordinates $y_i$ of $L''$ are bilinear functions of $ L, L'$.
\end{lem}
\Proof
$e_o \wedge x$ is not contained in $\pi$, hence does not belong to the line $\Ga = L * L'$,
and the first assertion is proven.

Write $ L'' = \sum_{i=1}^3 y_ i e_0 \wedge e_i = e_0 \wedge y$: then $L'' = e_o \wedge x$
if and only if it  contains $x$, or equivalently 
if and only if $L''$ is coplanar with $L$ and with $L'$, i.e. we have:
$$y_1 L_{23} - y_2  L_{13} + y_3  L_{12} = 0,$$
$$y_1 L'_{23} - y_2  L'_{13} + y_3  L'_{12} = 0.$$

The second assertion follows then from Cramer's rule,
$$ y_1 =  L_{13} L'_{12} - L'_{13} L_{12},  y_2 = L_{23}  L'_{12} - L'_{23} L_{12},$$
$$ y_3 =  L_{13}  L'_{12} - L'_{13} L_{12}.$$ 

\qed

We can now apply the lemma for the lines 
$ L \in Z, L' : = \nabla F (L) $, obtain a third line $L''$ which together with
$ L, L'$ yields a basis of $\PP^2_x$, under the assumption that
$F$ satisfies the weak Cayley equation, i.e., is a Cayley form.

Then, since the line $\Ga = L * L'$ is contained in $Z$, automatically we
obtain
$$D^2 F(L) (L, L) =  D^2 F(L) (L, L')= D^2 F(L) (L', L') = 0.$$

Hence follows immediately the following 

\begin{theo}\label{honest}
Let $F$ be a Cayley form. Then $F$ is a honest Cayley form if moreover
for each $ L \in Z$ the following equations hold:
$$D^2 F(L) (L'', L) =  D^2 F(L) (L'', L')= D^2 F(L) (L'', L'') = 0.$$

I.e., if and only if the above three polynomials,
whose coefficients have degree 2 or 3 in the coefficients of $F$, belong to 
the ideal $(Q, F)$ of $Z$.
\end{theo}
\Proof
The entries of the matrix $D^2 F(L)$ are linear in the coefficients of $F$,
as well as the coordinates of $L'$, while the coordinates of $L$
are homogeneous of degree 0 in the coefficients of $F$.
Since the Pl\"ucker coordinates $y_i$ of $L''$ are bilinear functions of $ L, L'$,
they are linear in the coefficients of $F$.

Hence the  three  equations are homogeneous in the coefficients of $F$,
of respective degrees $2,3,3  $. 

\qed

The next natural question is whether we can obtain from the above  theorem equations 
which hold $ mod (Q)$: we show that the answer is negative, already in the example
of a chain of three lines.

In this case, as we observed,   a Cayley form is
$$ F :=  p_{01} p_{02} p_{23},$$
and $F_2$ is here unique, equal to
$$   F_2 = F - \frac{1}{3} p_{02}  Q .$$

We set $L : =  p$, hence $L'  = \nabla F  - \frac{1}{3} p_{02} \nabla  Q - \frac{1}{3}  Q \nabla p_{02}   $,
and the equations determining $L''$ are 
$$   \L L'',  L \R  = 0, 0 =  \L L'',   \nabla F - \frac{1}{3}  Q \nabla p_{02} \R=y_1 p_{02}  p_{23}+ y_2 p_{01}  p_{23} - \frac{1}{3} Q y_2 . $$

These yield (modulo $Q$)
$$ y_1 p_{02} +  y_2 p_{01} = 0,\ \  y_3 p_{12} +   (p_{01} p_{23} +  p_{02} p_{13}) = 0,$$
hence as solution (modulo $Q$)

$$ y_1 =  p_{01} p_{12} , \ \  y_2 =  - p_{02} p_{12} ,\ \  y_3 =  - (p_{01} p_{23} + p_{02} p_{13}).$$
Observe now that, denoting  by $Q (q, q')$ the bilinear form associated to $Q$,
namely, $Q (q, q') : = \L q, q'  \R$, we have  $ Q (L'', L'') \equiv 0 $ and also  $Q(L,L'') \equiv Q (L', L'') \equiv 0 \  ( \ mod Q)$. Further $ Q (L, L') \equiv 0 $ on $Q$ (and also $ Q (L', L') \equiv 0 $ since we use $F_2$ for defining $L'$) while 
$ Q(L,L) \equiv 0$ holds tautologically on $Q$.

Since we are considering a point  $L = p \in Q$, when we look at the equations $D^2 F_2(L) (L'', L) =  D^2 F_2(L) (L'', L')= D^2 F_2(L) (L'', L'') = 0$,
we may replace it by  the simpler equations $D^2 F(L) (L'', L) =  D^2 F(L) (L'', L')= D^2 F(L) (L'', L'') = 0.$

Because $$ D^2 ( p_{02} Q) (q, q') =  2\  p_{02} Q (q, q') + q_{02} Q (p, q') + q'_{02} Q (p, q) .$$

 Now, whereas $$\frac{1}{2} D^2 F(L) (L'', L) =  p_{01} [y_2  p_{23}] + p_{02} [ y_1 p_{23}] + p_{23} [ y_2 p_{01} + y_1 p_{02}]  \equiv 0, $$
$$\frac{1}{2}  D^2 F(L) (L'', L'') =  p_{23} [ y_1 y_2] = - p_{12} ^2 p_{23} p_{01}  p_{02} = - p_{12} ^2  F  $$

which is not identically zero modulo $Q$. We have therefore shown

\begin{prop} Consider the equations in theorem \ref{honest}
for a honest Cayley form:
$$D^2 F(L) (L'', L) =  D^2 F(L) (L'', L')= D^2 F(L) (L'', L'') = 0.$$

If we take a chain $C$ of three lines in $\PP^3$, then the representative $F_2$ is unique, and for any choice of a Cayley form 
for $C$ these equations belong to 
the ideal $(Q, F)$ of $Z$, but not to the ideal of $Q$.
\end{prop}
\bigskip

{\bf Acknowledgements.}

I would like to thank Mark Green for stimulating email correspondence at the 
beginning of this research.

The research project was  begun during a visit of Hubert Flenner to Pisa, supported 
by an EEC research contract, 
and continued at 
the Max Planck Institut, Bonn, in July 1998, during a special activity dedicated to the
 memory of Boris Moishezon.

I would like to  thank  Hubert for all the efforts he devoted to the enterprise of deciding 
whether Chow varieties are defined by equations of
degree 2 and 3; at a certain point he however  decided that this project was too
 difficult to be pursued in its greatest generality
and stepped out. 

I would also like to thank  Michel Brion for interesting discussions on  representation 
theoretic aspects of the higher dimensional case,
and Igor Dolgachev for useful comments and for spotting a computational mistake in corollary \ref{examples} .


\end{document}